\documentclass[a4paper, 11pt]{amsart}
\usepackage{ psfrag, epsfig,amsmath,amssymb,amsthm,array,multicol,verbatim,moreverb,mathrsfs}
\usepackage{chicago}
\usepackage{setspace}

\newtheorem{theorem}{Theorem}

\newtheorem{lemma}{Lemma}
\newtheorem{corollary}{Corollary}
\newtheorem{assumption}{Assumption}
\newtheorem{proposition}{Proposition}
\newtheorem{ex}{Example}
\newenvironment{example}{\begin{ex}\rm}{\smallskip\end{ex}}
\newtheorem{rem}{Remark}
\newenvironment{remark}{\begin{rem}\rm}{\smallskip\end{rem}}
\def\proof#1{{\bf #1.}}
\def\endproof{\hfill$\Box$}

\input cyracc.def 

\def\P{\mathcal P}

\newcommand{\bb}{\mathbb}

\def \R {{\bb R}}
\def \real {\R}

\def \Z {{\bb Z}}

\def \ot {\otimes}
\def \op {\oplus}

\def \Xcal{{\mathcal X}}
\def \Kcal{{\mathcal K}}
\def \Gcal{{\mathcal G}}
\def \Vcal{{\mathcal V}}

\def \Ecal{{\mathcal E}}
\def \Ccal{{\mathcal C}}

\def\ind{{\rm 1\hspace{-0.90ex}1}}
\def\prob{\mathbb P}
\def\esp{\mathbb E}

\def\1{{\bf 1}}
\def\0{{\bf 0}}

\title[Tail Asymptotics for Discrete Event Systems]{Tail Asymptotics for Discrete Event Systems}
\author[Marc Lelarge]{Marc Lelarge}

\begin{document}

\maketitle 


\begin{abstract}
In the context of communication networks, the framework of
stochastic event graphs allows a modeling of control mechanisms
induced by the communication protocol and an analysis of its
performances. We concentrate on the logarithmic tail asymptotics of
the stationary response time for a class of networks that admit a
representation as (max,plus)-linear systems in a random medium. We
are able to derive analytic results when the distribution of the
holding times are light-tailed. We show that the lack of
independence may lead in dimension bigger than one to non-trivial
effects in the asymptotics of the sojourn time. We also study in
detail a simple queueing network with multipath routing.
\end{abstract}

\section{Introduction}

In this paper, we study tail asymptotics of the form
\begin{eqnarray*}
\lim_{x\rightarrow \infty}\frac{1}{x}\log \prob(Z>x)=-\theta^*,
\end{eqnarray*}
where the random variable $Z$ corresponds to a "global" state
variable associated to a (max,plus)-linear system. We only deal with
light-tailed distributions, i.e. distribution functions that are
decaying exponentially fast. The simplest example of random variable
covered by our results is the stationary waiting time in a FIFO
$GI/GI/1$ queue. This case has been extensively studied in the
literature and much finer estimates are available, see the
complementary works \cite{igl} and \cite{pak}.

In recent years, there has been some interest in extending this
result to networks of queues. \cite{csc:intree} considers an
intree network and uses large deviations techniques to show that the
queue length distributions have an exponentially decaying tail. \cite{ganan} obtains the decay rate of the tail
distribution for two exponential server queues in series fed by
renewal arrivals. In
\cite{bpt}, the decay rate of the stationary waiting time and queue
length distributions at each node is computed in an acyclic network in the
context of quite general arrival and service processes. Literature
on large deviations of queueing networks with feedback is rare and
confined to the setting of networks described by finite-dimensional
Markov processes, see \cite{dupel}, \cite{dupelwei} and the recent works \cite{ign}, \cite{ignaop}. Moreover, these works concentrate on
local large deviations and cannot handle the large deviations of the
network in its stationary regime. The large deviations asymptotics
of queueing systems are difficult to analyze because they are
dynamical systems with discontinuities. To the best of our
knowledge, there is no rigorous result on the large deviations of
non-exponential networks with feedback in their stationary regime.

In this paper, we consider a class of networks that admit a
(max,plus)-linear representation. This class contains the stochastic
event graphs (which can be used to model window-based congestion
control mechanism like TCP) and hence our results give the tail
asymptotics of the steady state end-to-end response times of these
networks. We should stress that the results of this paper are not
restricted to this sub-class and we give an example of a network
with multipath routing that is covered by our framework.

From a mathematical point of view, we study $Z$ the stationary
solution of a (max,plus)-linear recursion. Precise results
concerning large deviations of products of random topical operators
have been obtained in \cite{too}. However very restrictive
conditions are required on the coefficients of the matrix and only
the irreducible case is studied in \cite{too}. Here we do not assume
these requirements to be fulfilled and we show that under mild
assumptions on the matrix structure, the tail behavior of $Z$ is
explicitly given and can be computed (or approximated) in practical
cases.

In the next section, we first give the general (max,plus) framework,
with some examples of queueing networks. Then we give the stochastic
assumptions and the tail asymptotics of the stationary solution of
the (max,plus)-linear recursion is derived in Theorem \ref{the:main}
which clearly extends the case of the single server queue. Theorem
\ref{the:2} gives a more explicit form of the rate of exponential
decay.

In Section \ref{sec:queue}, we study two queueing applications.
First we consider a system of two queues in tandem and show that
when the service times at both queues are identical, then depending
on the value of the intensity of the arrival process there is a
phase transition in the behavior of the network reaching a large
end-to-end delay (Proposition \ref{prop:tand}). Then we study in
detail a simple example of queueing networks with resequencing.
Multipath routing has recently received some attention in the
context of both wired and wireless communication networks. By
sending data packets along different paths, multipath routing can
potentially help balance the traffic load and reduce congestion
levels in the network, thereby resulting in lower end-to-end delay.
We show how our framework can model such mechanisms and give
analytical insights.

Sections \ref{sec:pmain}, \ref{sec:p2} contain the proofs of Theorem
\ref{the:main} and Theorem \ref{the:2} respectively.
We give some further directions of research in the conclusion.

\section{Tails for discrete event systems}

In this paper we consider open systems with a single input marked
point process $N=\{(T_n, A_n, B_n)\}_{-\infty<n<\infty}$, where in a
queueing context the sequence $\{T_n\}$ describes the arrival times
of customer in the network and $\{A_n,B_n\}$ carries the information
related to the $n$-th customer (like its service time at the
different stations, routing decisions). We give a precise
description of the dynamic of the system in the next section and of
the random variable $Z$ for which we derive the logarithmic tail
asymptotics.

\subsection{(Max, plus)-linear systems}
\label{sec:maxpluslinear}

The (max, plus) semi-ring $\real_{\max}$ is the set $\real\cup
\{-\infty\}$, equipped with $\max$, written additively (i.e.,
$a\oplus b=\max(a,b)$) and the usual sum, written multiplicatively
(i.e., $a\otimes b = a+b$). The zero element is $-\infty$. For
matrices of appropriate sizes, we define $(A\oplus
B)^{(i,j)}=A^{(i,j)}\oplus B^{(i,j)}:=\max(A^{(i,j)},B^{(i,j)})$ and
$(A\otimes B)^{(i,j)} = \bigoplus_{k}A^{(i,k)}\otimes
B^{(k,j)}:=\max_{k}(A^{(i,k)}+B^{(k,j)})$. By convention if $A$ is a
matrix and $c\in \real_{\max}$, then $(A\ot c)^{(i,j)}:=
A^{(i,j)}\ot c$.

Let $s$ be an arbitrary fixed natural number. We assume that we are
given with a sequence of matrices with non-negative coefficients:
$A_n$ of size $s\times s$ and $B_n$ of size $s\times 1$. To the
sequences $\{A_n\}_n$, $\{B_n\}_n$, and $\{T_n\}_n$, we associate
the following (max, plus)-linear recurrence:
\begin{eqnarray}
\label{eq:rec} \Xcal_{n+1}=A_{n+1}\ot \Xcal_n\op B_{n+1}\ot T_{n+1},
\end{eqnarray}
where $\{\Xcal_n, n\in \Z\}$ is a sequence of state variables of
dimension $s$. In Examples \ref{exa:ss1}, \ref{exa:ta1},
\ref{exa:fj1}, we derive the explicit form of this recursion for the
single server queue, queues in tandem and a fork join system. We
refer to these examples to get an interpretation of the various
quantities.

The stationary solution of this equation is constructed as follows.
We write
\begin{eqnarray}
\label{def:Y} Y_{[m,n]} := \bigoplus_{m\leq k \leq n}
D_{[k+1,n]}\otimes B_{k} \ot T_k, 
\end{eqnarray}
where for $k<n$, $D_{[k+1,n]}= \bigotimes_{j = n}^{k+1}A_j=A_{n}
\otimes \dots \otimes A_{k+1}$ and $D_{[n+1,n]}=E$, the identity
matrix (the matrix with all its diagonal elements equal to 0 and all
its non-diagonal elements equal to $-\infty$). It is readily checked
that $Y_{[m,m]}=B_m\otimes T_m$, and for all $n\geq m$,
\begin{eqnarray*}
Y_{[m,n+1]}=A_{n+1}\otimes Y_{[m,n]} \oplus B_{n+1}\otimes T_{n+1}.
\end{eqnarray*}

In view of (\ref{def:Y}), the sequence $\{Y_{[-n,0]}\}$ is
non-decreasing in $n$, so that we can define the stationary solution
of (\ref{eq:rec}),
\begin{eqnarray*}
Y_{(-\infty,0]} := \lim_{n\rightarrow \infty}Y_{[-n,0]}\leq \infty.
\end{eqnarray*}
We define the {\em stationary maximal dater} by
\begin{eqnarray}
\label{def:maxdat} 0\leq Z:=Z_{(-\infty,0]} =\bigoplus_{1\leq i\leq
s}Y_{(-\infty,0]}^{(i)}-T_0\leq \infty.
\end{eqnarray}
The following writing for the stationary maximal dater shows the
similitude with the traditional stationary workload of a single
server queue:
\begin{eqnarray}
\label{eq:Zsup}Z =\sup_{n\leq 0} \{S_n - (T_0-T_n)\},
\end{eqnarray}
where the process $\{S_n\}_{n\leq 0}$ is defined by
\begin{eqnarray}
\label{eq:defS_n}S_n := \bigoplus_{1\leq i\leq s} 
\left(D_{[n+1,0]}\otimes B_n\right)^{(i)}.
\end{eqnarray}

\begin{example}\label{exa:ss1}
Consider a FIFO single server queue where $T_n$ is the arrival time
of the $n$-th customer and $\sigma_n$ is its service time. Equation
(\ref{eq:rec}) is then the standard Lindley's recursion,
\begin{eqnarray*}
\Xcal_{n+1}&=&\max(\Xcal_n +\sigma_{n+1}, T_{n+1}+\sigma_{n+1})\\
&=& \sigma_{n+1}\ot \Xcal_n\op \sigma_{n+1}\ot T_{n+1}.
\end{eqnarray*}
The interpretation of $\Xcal_n$ is the departure time of the $n$-th
customer from the queue. Note that in this case, we have $S_n =
\sum_{i=n}^0 \sigma_i$ and $Z$ is the stationary workload.
\end{example}

\begin{example}\label{exa:ta1} Consider now a system of two queues
in tandem, where $T_n$ is the arrival time of the $n$-th customer in
the system and $\sigma^{(i)}_n$ is its service time at queue $i$,
for $i=1,2$. Then Equation (\ref{eq:rec}) is given by
\begin{eqnarray*}
\left(\begin{array}{c} \Xcal_{n+1}^{(1)}\\
\Xcal_{n+1}^{(2)}\end{array}
\right)&=&\left(\begin{array}{cc} \sigma_{n+1}^{(1)}&-\infty\\
\sigma_{n+1}^{(1\ot 2)}& \sigma_{n+1}^{(2)}\end{array}
\right)\left(\begin{array}{c} \Xcal_{n}^{(1)}\\
\Xcal_{n}^{(2)}\end{array} \right)\oplus \left(\begin{array}{c} \sigma_{n+1}^{(1)}\\
\sigma_{n+1}^{(1\ot 2)}\end{array} \right)T_{n+1},
\end{eqnarray*}
where we used the shorthand notation $\sigma_{n+1}^{(1\ot
2)}=\sigma_{n+1}^{(1)}\ot\sigma_{n+1}^{(2)}$. In this case
$\Xcal^{(1)}_n$ is the departure time of the $n$-th customer from
the first queue and $\Xcal^{(2)}_n$ from the second queue. Hence $Z$
is the stationary end-to-end delay of the network. Note that in this
case, we have
\begin{eqnarray}
\label{eq:Stand}S_n =\sup_{n\leq \ell\leq
0}\sum_{i=n}^\ell\sigma^{(1)}_i+\sum_{j=\ell}^0 \sigma^{(2)}_j.
\end{eqnarray}
\end{example}

\begin{example}\label{exa:fj1}
Consider the standard fork and join system as depicted (with Petri
net formalism) in Figure \ref{fig:forkjoin}. In this model, each
time a packet (say $k$) finishes its service $\sigma^{(1)}_k$ in
node $1$, there is one packet sent up and one packet sent down
simultaneously. The `up'-packet (`down'-packet) is then also the
$k$-th packet for node 2 (for node 3 respectively). The $k$-th
packet joins the queue of node 4 once both packets have left node 2
and 3 respectively. Each node is a standard $\cdot/G/1/\infty$
queue.
\begin{figure}[hbt]
\center \psfrag{t1}{\Huge{ $\sigma^{(1)}$}} \psfrag{t2}{\Huge{
$\sigma^{(2)}$}} \psfrag{t3}{\Huge{ $\sigma^{(3)}$}}
\psfrag{ts}{}
\resizebox{50mm}{!}{\includegraphics*{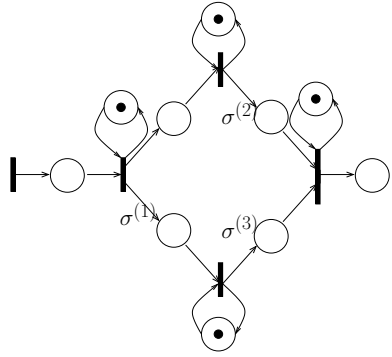}} \caption{Fork
and join model.} \label{fig:forkjoin}
\end{figure}

Let $\Xcal^{(i)}_n$ denotes the departure time of the $n$-th packet
from node $i$. We have the following equations:
\begin{eqnarray*}
\Xcal^{(1)}_{n+1} &=& (T_{n+1}\oplus \Xcal^{(1)}_n)\otimes \sigma^{(1)}_{n+1},\\
\Xcal^{(2)}_{n+1} &=& (\Xcal^{(1)}_{n+1}\oplus \Xcal^{(2)}_n)\otimes \sigma^{(2)}_{n+1},\\
\Xcal^{(3)}_{n+1} &=& (\Xcal^{(1)}_{n+1}\oplus \Xcal^{(3)}_n)\otimes \sigma^{(3)}_{n+1},\\
\Xcal^{(4)}_{n+1} &=& (\Xcal^{(2)}_{n+1}\oplus
\Xcal^{(3)}_{n+1}\oplus \Xcal^{(4)}_n)\otimes 0.
\end{eqnarray*}
This system is linear in the (max, plus) semi-ring $\real_{\max}$,
and we can write the recursion (\ref{eq:rec}) with the following
matrices:
\begin{eqnarray*}
A_n=\left(\begin{array}{cccc}
{\sigma^{(1)}_n}&-\infty&-\infty&-\infty\\
\sigma^{(1\ot 2)}_n&{\sigma^{(2)}_n}&-\infty&-\infty\\
\sigma^{(1\ot 3)}_n&-\infty&{\sigma^{(3)}_n}&-\infty\\
\sigma^{(1 \ot 2\oplus 3)}_n&\sigma^{(2)}_n&\sigma^{(3)}_n&{0}
\end{array}\right),&&
B_n=\left(\begin{array}{c}\sigma^{(1)}_n\\
\sigma^{(1\ot 2)}_n\\
\sigma^{(1\ot 3)}_n\\
\sigma^{(1\ot 2\oplus 3)}_n
\end{array} \right),
\end{eqnarray*}
where we used the shorthand notations, $\sigma^{(i\ot
j)}_n=\sigma^{(i)}_n\ot \sigma^{(j)}_n$ and $\sigma^{(i \ot j\oplus
k)}_n=\sigma^{(i)}_n\ot \left(\sigma^{(j)}_n\op
\sigma^{(k)}_n\right)$. In this case $Z$ is the stationary
end-to-end delay of the network.
\end{example}

We refer to \cite{bcoq} for other examples of (max,plus)-linear
networks (see also \cite{these} for an example showing how to model
window control mechanism).

\subsection{Tail asymptotics for the stationary solution}

We first need to give some conditions that ensure the stability of
the system, i.e. that the limit (\ref{def:maxdat}) is finite. Hence
we introduce the following assumption in order to apply first-order
Theorems of Section 7.3 of \cite{bcoq}:

\renewcommand{\theassumption}{{\bf (ST)}}
\begin{assumption}\label{ass:st}
{\em (Structure of $(A_n, B_n)$)}

\noindent The random variables $A_n,B_n$ have a fixed structure,
i.e. each entry of $A_n$ or $B_n$ is either a.s. equal to $-\infty$
or non-negative for all $n$. And each diagonal entry of $A_n$ is
non-negative.
\end{assumption}

\renewcommand{\theassumption}{{\bf (IA)}}
\begin{assumption}\label{ass:ia}
  {\em (independence assumption)}

\noindent We suppose that the sequences $\{(A_n,B_n)\}_n$ and
$\{\tau_n := T_{n+1}-T_n\}_n$ are mutually independent and each of
them consists of i.i.d.~random variables with finite means.
\end{assumption}

This assumption implies a law of large number for $\{S_{-n}\}$ defined
in (\ref{eq:defS_n}), namely,
\begin{eqnarray}\label{eq:stab}
\frac{S_{-n}}{n}\rightarrow_{n\to\infty} \gamma \quad \mbox{both
a.s.~and in $L_1$,}
\end{eqnarray}
where $\gamma$ is a constant referred to as the top Lyapunov
exponent of the sequence $\{A_n\}$ see Theorems 7.27 and 7.36 in
\cite{bcoq}.

\renewcommand{\theassumption}{{\bf (S)}}
\begin{assumption}\label{ass:s}
{\em (stability)}

\noindent We assume that $\gamma <\esp[\tau_1]=a$.
\end{assumption}

We have that under \ref{ass:ia} and \ref{ass:s} the maximal dater
$Z$ defined in (\ref{def:maxdat}) is almost surely finite.

We denote by $\0$ the vector with all its entries equal to $0$.
\renewcommand{\theassumption}{{\bf (SP)}}
\begin{assumption}\label{ass:sp}
{\em (Separability)}

\noindent We assume that we have for all $n$,
\begin{eqnarray*}
A_n \otimes\0 = B_n\oplus \0.
\end{eqnarray*}
\end{assumption}

This assumption ensures that for a solution $\Xcal_n$ of
(\ref{eq:rec}) with any initial condition: if $\Xcal_n\leq\0\otimes
T_{n+1}$ then the process $\Xcal_{n+1},\Xcal_{n+2},\dots$ does not
depend on the past $\Xcal_{n},\Xcal_{n-1},\dots$ Note that this
assumption is clearly satisfied in the examples described above
since we have $A_n \otimes\0 = B_n$. We refer to Section 2.2.4 of
\cite{these} for an example of network with $A_n \otimes\0 =
B_n\oplus \0\neq B_n$. In fact, Propositions 3 and 4 of \cite{these}
show that any FIFO event graph with a single input fits into our
framework. However this condition allows also to deal with some type
of networks with (random) routing as described in Section
\ref{sec:reseq}. This property of separability can be made precise
in a larger framework than (max,plus)-linear networks: this is the
class of monotone separable networks introduced in
\cite{bacfos:sat}.

We now give the stochastic assumptions that ensure that the random
variable $Z$ is light tailed:
\renewcommand{\theassumption}{{\bf (LT)}}
\begin{assumption}\label{ass:lt}
{\em (Light-tailed)}

\noindent Let
\begin{eqnarray*}
\eta = \sup\left\{\theta>0,\: \bigoplus_{i}\esp\left[e^{\theta
B^{(i)}_1}\right]<\infty\right\}.
\end{eqnarray*}
We assume that $\eta>0$.
\end{assumption}

We will always assume that Assumptions {\bf (ST), (IA), (S), (SP),
(LT)} hold. We are now in position to state our main result.
\begin{theorem}\label{the:main}
The following limit exists as an extended real number:
\begin{eqnarray}
\label{def:lambda}\Lambda_S(\theta) = \lim_{n\to
\infty}\frac{1}{n}\log \esp\left[ e^{\theta S_{-n}}\right].
\end{eqnarray}
We have
\begin{eqnarray}
\label{def:theta}\theta^* = \sup\{\theta>0,\:
\Lambda_S(\theta)+\Lambda_T(-\theta)<0\}>0,
\end{eqnarray}
where $\Lambda_T(\theta)=\log\esp\left[e^{\theta \tau_1}\right]$ and
the tail asymptotics of $Z$ is given by,
\begin{eqnarray*}
\lim_{x\to \infty}\frac{1}{x}\log\prob(Z>x) = -\theta^*.
\end{eqnarray*}
\end{theorem}

In the case of the single server queue, we have clearly
$\Lambda_S(\theta) = \log\esp[\exp(\theta \sigma_1)]$ and Theorem
\ref{the:main} is a standard result of queueing theory that goes
back to the work \cite{cramer} and in a queueing context
to \cite{igl}. We will give more comments on this theorem
in Section \ref{sec:queue}.

\subsection{More detailed results}\label{sec:det}

In this section we give a more explicit form for $\theta^*$. Without
loss of generality, we may assume that the matrices $A_n$ have the
following block structure:
\begin{eqnarray*}
\left(\begin{array}{ccccccc}
A_n(1,1)&|&-\infty&|&-\infty&|&-\infty\\
-&-&-&-&-&-&-\\
A_n(2,1)&|&A_n(2,2)&|&-\infty&|&-\infty\\
-&-&-&-&-&-&-\\
&\vdots&&\vdots&&\vdots&\\
-&-&-&-&-&-&-\\
A_n(d,1)&|&A_n(d,2)&|&&|&A_n(d,d)\\
\end{array}\right),
\end{eqnarray*}
where each $A_n(\ell,\ell)$ is an irreducible matrix.

\begin{theorem}\label{the:2}
Associated to the irreducible matrices $\{A_n(\ell,\ell)\}$, we
define the following function:
\begin{eqnarray*}
\Lambda_\ell(\theta) = \lim_{n\rightarrow \infty} \frac{1}{n}
\log\esp\left[ e^{\theta (A_n(\ell,\ell)\otimes \dots \otimes
A_1(\ell,\ell))^{(u,v)}}\right],
\end{eqnarray*}
where the limit exists in $\real\cup\{\infty\}$ and is independent
of $u,v$. Then we have $\theta^*= \min\{\eta,\theta^\ell\}$ where
the $\theta^\ell$'s are defined as follows
\begin{eqnarray*}
\theta^\ell = \sup\{\theta>0,\:
\Lambda_\ell(\theta)+\Lambda_T(-\theta)<0 \}.
\end{eqnarray*}
\end{theorem}

In the case of a single server queue with exponentially distributed
service times, we have $\eta>\theta^*$ and this property remains
valid for a large class of distributions. However, we show in the
next section that as soon as we consider a network (i.e. with at
least 2 nodes) then the parameter $\eta$ can play a role even with
exponentially distributed service times.

We first give a framework where $\eta$ cannot play any role. 
Given a vector $v=(v^{(1)},\ldots, v^{(K)})$, we call a (max, plus)
expression $\P$ a polynomial in $v$ of unit maximum degree if it has
the form
\[
\P =\bigoplus_j \bigotimes_{k\in \Kcal_j} v^{(k)},
\]
where $\Kcal_j\subset [1,K]$.
\begin{corollary}\label{cor}
If there exists a sequence of random variables
$\{\sigma_n=(\sigma^{(1)}_n,\dots,\sigma^{(K)}_n)\}_n$ such that
\begin{enumerate}
\item the components of $\sigma_n$ are independent of each other;
\item for all $i$, there exists $k$ such that $A_n^{(k,k)}=\sigma^{(i)}_n$;
\item each entry of $A_n$ that is not $0$ or $-\infty$ is a
polynomial (in $\real_{\max}$) in $\sigma_n$ of unit maximal degree.
\end{enumerate} Then we have $\theta^*= \min\{\theta^\ell\}$.
\end{corollary}

In a queueing context, the sequence of matrices $\{A_n(\ell,\ell)\}$
corresponds to a specific "component" of the network. It is
well-known that the stability of such a network is constraint by the
"slowest" component. Here we see that in a large deviations regime,
if each component is independent of each other, then the "bad"
behavior of the network is due to a "bottleneck" component (which is
not necessarily the same as the "slowest" component in average).

\begin{remark}\label{rem:bruno}
In the framework of last Corollary, the tail asymptotics
for $Z$ under heavy-tailed (more precisely subexponential)
assumptions (i.e. when Assumption {\bf (LT)} is not satisfied) has
been derived in \cite{bfl}. In this case, the exact asymptotics
(i.e. not in the logarithmic scale) are derived and the Lyapunov
exponents of the sub-matrices $\{A_n(\ell,\ell)\}$ appear. These
exponents are known to be hard to compute \cite{bgt}. Similarly in
the light-tailed case, we see that the asymptotics in the
logarithmic scale is given by the functions $\Lambda_\ell(\theta)$.
The computation of these functions is not easy, in particular when
the network has some feedback mechanism. More formally, we will see
that this function is convex and its right-derivative at zero is
exactly the Lyapunov exponent of the sub-matrices
$\{A_n(\ell,\ell)\}$, so that knowing the function
$\Lambda_\ell(\theta)$ allows to determine the Lyapunov exponent.
In particular, at the level of generality considered in this paper, our result cannot be made in a more explicit form. We will
see in the next section several examples for which the value of
$\theta^*$ has a simple expression in term of the parameters of the
problem. It is
interesting to note that, as opposed to the
heavy-tailed case, exact tail asymptotics for $Z$ under light-tailed
assumptions seem  to be out of reach in the general framework of
(max,plus)-linear networks.
\end{remark}

\begin{example}
Going back to the fork and join system described in Example
\ref{exa:fj1}, we see that the irreducible matrices are
one-dimensional and we have for $\ell=1,2,3$,
\begin{eqnarray*}
\Lambda_\ell (\theta) = \log \esp\left[ e^{\theta
\sigma^{(\ell)}_1}\right].
\end{eqnarray*}
Hence if $\theta^\ell$ denotes the exponential rate of decay for the
single server queue fed by $\{T_n, \sigma^{(\ell)}_n\}$, then we
have $\theta^*=\min\{\theta^\ell\}$ in the case where each sequence
of service times at each station are independent of each other.
\end{example}

\section{Queueing Applications}\label{sec:queue}

\subsection{The impact of dependence}

In view of (\ref{eq:Zsup}), $Z$ is the supremum of a random process
with negative drift and to make the connection with the existing
literature, we state the following result (for a proof we refer to
\cite{tamsn} Corollary 3.2):
\begin{proposition}\label{cor}
Under Assumptions {\bf (IA)} and {\bf (S)} and if
\begin{enumerate}
\item the sequence $\{S_{-n}/n\}$ satisfies a large deviation principle (LDP) with a good rate
function I;
\item there exists $\epsilon>0$ such that
$\Lambda_S(\theta^*+\epsilon)<\infty$,
\end{enumerate}
where $\theta^*$ is defined as in (\ref{def:theta}). Then we have
\begin{eqnarray}
\label{ssq}\lim_{x\rightarrow \infty}\frac{1}{x}\log \prob(Z>x) =
-\theta^* = -\inf_{\alpha>0} \frac{I(\alpha)}{\alpha}.
\end{eqnarray}
\end{proposition}
This kind of result has been extensively studied in the queueing
literature and follows directly  from the work \cite{duffy}.
However, we see that considering the moment generating function
instead of the rate function allows us to get a more general result
than (\ref{ssq}) since we do not require the assumption on the tail
(which is essential for (\ref{ssq}) to hold see \cite{duffy}).
Indeed this assumption ensures that the tail asymptotics of
$\prob(S_n-(T_0-T_n)>nc)$ for a single $n$ value cannot dominate
those of $\prob(M>x)$. In this case, equation (\ref{ssq}) has a nice
interpretation: the natural drift of the process $S_n-(T_0-T_n)$ is
$(\gamma-a) n$, where $\gamma-a<0$. The quantity $I(\alpha)$ can be
seen as the cost for changing the drift of this process to
$\alpha>0$. Now in order to reach level $x$, this drift has to last
for a time $x/\alpha$. Hence the total cost for reaching level $x$
with drift $\alpha$ is $xI(\alpha)/\alpha$ and the process naturally
chooses the drift with the minimal associated cost. This can be made
precise in some cases by a conditional limit theorem that
characterizes the most likely path.

In this section we show that this interpretation might be misleading
in a queueing context. We consider a very simple example: a system
of two queues in tandem. We assume that the sequence
$\{(\sigma^{(1)}_n,\sigma^{(2)}_n)\}_n$ is a sequence of i.i.d.
random variables with
$\max\{\esp[\sigma^{(1)}_1],\esp[\sigma^{(2)}_1]\}<a$ and
$\esp[\exp\theta( \sigma^{(1)}_1+\sigma^{(2)}_1)]$ finite in a
neighborhood of the origin.

If the service times of station 1 and station 2 are independent of
each other, then the most likely cause of a given customer suffering
a large delay is that a large number of its immediate predecessors
require service times in excess of their inter-arrival times at one
of the station. However in the case where the service times are the
same at both stations, we show that depending on the intensity of
the arrival process $\lambda$, two situations may occur:
\begin{enumerate}
\item if $\lambda<\lambda_c$, then the most likely reason that a
given customer suffers a large delay is that its own service time is
large;
\item if $\lambda>\lambda_c$, then the tail asymptotic of the
end-to-end delay is the same as in the independent case.
\end{enumerate}
Let consider first the case where $\sigma^{(1)}_n$ and
$\sigma^{(2)}_n$ are independent. We are in the framework of
Corollary \ref{cor}. Hence if we denote by $\theta ^{(i)}$ the
exponential rate of decay for the tail asymptotics of the stationary
workload of a single server queue with arrival times $T_n$ and
service times $\sigma^{(i)}_n$, then we have
\begin{eqnarray*}
\lim_{x\to\infty} \frac{1}{x} \log \prob (Z>x) =-\min (\theta^{(1)},
\theta^{(2)}).
\end{eqnarray*}
This result has been obtained in \cite{ga}. 
In words, we can say that the large deviation of the end-to-end
delay in a system of two queues in tandem with independent service
times is dominated by the "worst" one.

Consider now the case where $\sigma^{(1)}_n=\sigma^{(2)}_n$ for all
$n$ and the sequence $\{\sigma^{(1)}_n\}_n$ is a sequence of i.i.d.
random variables exponentially distributed with mean $1/\mu$. We
assume also that the arrival process is Poisson with rate
$\lambda<\mu$. Then a direct application of Theorem \ref{the:2}
gives,
\begin{proposition}\label{prop:tand}
In the previous framework, we have
\begin{eqnarray*}
\lim_{x\to\infty} \frac{1}{x} \log \prob (Z>x) =-\theta^*,
\end{eqnarray*}
with\begin{eqnarray*}
\lambda\leq \mu/2 &\Rightarrow & \theta^*= \mu/2,\\
\lambda>\mu/2 &\Rightarrow & \theta^*= \mu-\lambda.
\end{eqnarray*}
\end{proposition}
This proposition completes the result in \cite{ga}. For small values
of $\lambda$, the tail of the end-to-end delay is determined by the
total service requirement of a single customer whereas when
$\lambda>\mu/2$, it is the same as in the independent case.

This shows that the behavior of tandems differs from that of a
single server queue. In particular \cite{venkat} shows
that for $GI/GI/1$ queues, the build-up of large delays can happen
in one of two ways. If the service times have exponential tails,
then it involves a large number of customers (whose inter-arrival
and service times differ from their mean values). This behavior is
analogous to that of tandems where the service times are independent
at each station or if the intensity of the arrival process is
sufficiently large. If the service times do not have exponential
tails, then large delays are caused by the arrival of a single
customer with large service requirement. In contrast, we see that a
single customer can create large delays in tandems even under the
assumption of exponential service times, if the intensity of the
arrival is sufficiently low.

\subsection{A case of study: queueing network with
resequencing}\label{sec:reseq}

The aim of this section is to show that the results of this paper
are not restricted to the class of event graphs and that our
framework can deal with complex synchronization problems encountered
in applications.

In many distributed applications (e.g., remote computations,
database manipulations, or data transmission over a computer
network), information integrity requires that data exchanges between
different nodes of a system be performed in a specific order.
However, due to random delays over different paths in a system, the
packets or updates may arrive at the receiver in a different order
than their chronological order. In such a case, a buffer (with
infinite capacity) at the receiver has to store disordered packets
temporarily. There is an extensive literature on resequencing
problem and we refer the interested reader to the survey
\cite{bacmak} (see also the more recent work \cite{jmg}).

\begin{figure}[hbt]
\center \resizebox{80mm}{!}{\includegraphics*{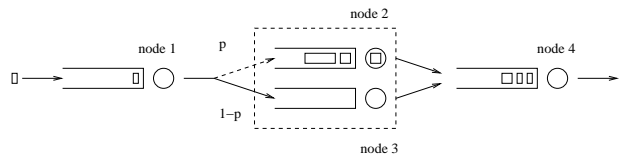}}
\caption{Resequencing problem.} \label{fig:reseq}
\end{figure}

We consider a simple queueing model of disordering, namely a set of
2 parallel single server queueing station $.|GI|1$ with renewal
arrivals under probabilistic state-independent routing. This model
constitutes an ersatz of the very complex situation one seek
investigate. While the details of any protocol have been eliminated,
the essence of network behavior (i.e. disordering) is preserved. In
the sequel, we shall thus consider the model described in Figure
\ref{fig:reseq}.
We assume that packets arrive in the first queue according to a
renewal process $\{T_n\}$. We will model the routing at node 1 by a
Bernoulli routing: with probability $p$ (resp. (1-p)) the packet is
sent up to node 2 (resp. down to node 3). Once packet $k$ reaches
the receiver, it leaves the system if all packets $j$ with $j<k$
have already left the system. Otherwise it stays in the resequencing
buffer, where it waits for the packets with number less than $k$.

This model is similar to the standard fork and join system described
in Example \ref{exa:fj1}. Even if the routing mechanism is not the
same here, we now show how we can adapt Example \ref{exa:fj1}. In
order to model the desired routing mechanism we will use the idea of
clones, i.e., packets that behave like real packets except that they
never require any service time: their service time is null. Suppose
that the real route of packet $k$ is up. Then at the end of its
service in the first node, a clone is sent to node 3. Since
$\sigma^{(3)}_k=0$, the departure time of packet $k$ from node 3 is
$\Xcal^{(3)}_{k}=\max (\Xcal^{(1)}_{k},\Xcal^{(3)}_{k-1})$.
Similarly, if the real route of packet $k$ is down, then a clone is
sent up. In both cases the ``real'' packet $k$ joins the queue of
node 4 once ``real'' packet $k-1$ has joined it (and not before). In
particular packets are ordered when they leave node 4.

Let $\{\zeta_n=(\zeta^{(1)}_n,\dots \zeta^{(3)}_n)\}_n$ be a
sequence of i.i.d. mutually independent random variables with finite
mean and such that $\bigoplus_i \esp\left[\exp(\theta
\zeta^{(i)}_1\right]<\infty$ for $\theta$ in a neighborhood of the
origin. Let $\{r_n\}_{n\in\Z}$ be a sequence of i.i.d. random
variables, independent of everything else, with values in $\{2,3\}$.
We write $\prob(r_n=2)= 1-\prob(r_n=3)=:p$, and assume that $0<p<1$.
In order to apply our idea of clones, we consider the
(max,plus)-recursion defined in Example \ref{exa:fj1} and we define
\begin{eqnarray*}
\sigma^{(1)}_n :=\zeta^{(1)}_n,\: \sigma^{(2)}_n
:=\zeta^{(2)}_n\ind_{\{r_n=2\}},\:
\sigma^{(3)}_n :=\zeta^{(3)}_n\ind_{\{r_n=3\}},\: \sigma^{(4)}_n
:=0.
\end{eqnarray*}

We see that our system satisfies the assumptions of Theorem
\ref{the:2}. In this case $Z$ is the end-to-end sojourn time of a
packet (taking into account the resequencing delay). The following
proposition is a direct application of Theorem \ref{the:2}. In
particular, note that $\eta\geq \theta^\ell$, hence we have
\begin{proposition}
In the previous framework, we have
\begin{eqnarray*}
\lim_{x\rightarrow \infty} \frac{1}{x} \log \prob(Z>x) =
-\theta^*<0,
\end{eqnarray*}
where $\theta^*= \min\{\theta^\ell\}$ and the $\theta^\ell$'s are
defined as follows
\begin{eqnarray*}
\theta^\ell = \sup\{\theta>0,\:
\Lambda_\ell(\theta)+\Lambda_T(-\theta)<0 \},
\end{eqnarray*}
with
\begin{eqnarray*}
\Lambda_1(\theta ) &=& \log \esp\left[ e^{\theta
\zeta^{(1)}_1}\right],\\
\Lambda_2(\theta ) &=& \log \left(p\esp\left[ e^{\theta
\zeta^{(2)}_1}\right]+1-p\right),\\
\Lambda_3(\theta ) &=& \log \left((1-p)\esp\left[ e^{\theta
\zeta^{(3)}_1}\right]+p\right).
\end{eqnarray*}
\end{proposition}

Then it is possible to make some optimizations. Assume to simplify
that $\zeta^{(1)}_n=0$ for all $n$ and that the sequence
$\{\zeta^{(2)}_n\}$ and $\{\zeta^{(3)}_n\}$ are independent
sequences of i.i.d. random variables exponentially distributed with
respective mean $1/\mu_2$ and $1/\mu_3$. We assume also that the
arrival process is Poisson with rate $\lambda>\max(\mu_2,\mu_3)$,
i.e. we cannot send all packets to one node. In this case, we have
to load-balance the traffic in order to have a stable system. We
have to chose $p$ such that $\lambda<\min(\mu_2/p,\mu_3/(1-p)\}$. In
order to do so efficiently, we would like to maximize $\theta^*$.
Here, we have $\theta_2 = \mu_2-\lambda p$ and $\theta_3 =\mu_3
-\lambda(1-p)$, hence the optimal value of $p$ is given by
\begin{eqnarray*}
p =\frac{1}{2}\left( \frac{\mu_2-\mu_3}{\lambda}+1\right)\mbox{ and
then, } \theta^*=\frac{\mu_2+\mu_3-\lambda}{2}.
\end{eqnarray*}
In particular, in the symmetric case $\mu_1=\mu_2$, we find $p=1/2$
which is in accordance with standard results of resequencing
literature \cite{jmg}.

This model is certainly an oversimplification of the reality but it
is intended to be also of pedagogical interest. We should stress
that it could incorporate quite a few sophistication to enable to
take into account general distributions, more than 2 paths or window
control mechanism (where at any time, the number of packets in the
network is upper bounded by the window size). In this cases,
computations are much more complex and numerical approximations have
to be made.

\section{Proof of Theorem \ref{the:main}}\label{sec:pmain}

We first prove the existence of the moment generating function
$\Lambda_S$ given by the limit (\ref{def:lambda}). Then we prove
that $\theta^*$ defined by (\ref{def:theta}) is positive and then we
derive the tail asymptotics for $Z$.

Note that under our assumptions we have for $n\geq 0$,
\begin{eqnarray}\label{eq:invtime}
S_n =_d \bigoplus_{1\leq i\leq s} \left(D_{[1,n]}\otimes B_0\right)^{(i)},
\end{eqnarray}
where the equality is in distribution. Hence for the simplicity of
notation, we will "inverse time" and index the processes by
non-negative indexes. Hence $S_n$ is now given by the right-hand
term of (\ref{eq:invtime})
and more generally we define for $u\leq v$,
\begin{eqnarray*}
S_{[u,v]}=\bigoplus_{1\leq i\leq s} \left(D_{[u+1,v]}\otimes B_u\right)^{(i)}.
\end{eqnarray*}

\subsection{Computation of the moment generating function}


Note that thanks to Assumption {\bf (ST)}, we have with $i^*$ defined
by $S_n = \left(D_{[1,n]}\otimes B_0\right)^{(i^*)}$, then
\begin{eqnarray*}
S_{n+1}\geq A_{n+1}^{(i^*,i^*)}+S_n\geq S_n.
\end{eqnarray*}
In particular the process $\{S_n\}$ is non-decreasing and we have the
following rewriting:
\begin{eqnarray*}
S_{[u,v]}=\bigoplus_{1\leq i\leq s} \bigoplus_{u\leq k\leq
v}\left(D_{[k+1,v]}\otimes B_k\right)^{(i)}.
\end{eqnarray*}
The following lemma shows that the process has a subadditive property.
\begin{lemma}\label{lem:aux}
We have for $n,m\geq 0$,
\begin{eqnarray*}
S_{n+m}\leq S_n + S_{[n+1,n+m]}.
\end{eqnarray*}
\end{lemma}
\proof{Proof} We have by definition,
\begin{eqnarray*}
\bigoplus_{0\leq k\leq n}D_{[k+1,n]}\otimes B_k \leq \0\ot S_n,
\end{eqnarray*}
hence by monotonicity, we have
\begin{eqnarray*}
A_{n+1} \ot\bigoplus_{0\leq k\leq n}D_{[k+1,n]}\otimes B_k &\leq&
A_{n+1}\ot \0 \ot S_n,\\
\bigoplus_{0\leq k\leq n}D_{[k+1,n+1]}\otimes
B_k&\leq&\left(B_{n+1}\op \0\right) \ot S_n,
\end{eqnarray*}
iterating we get:
\begin{eqnarray*}
{\bigoplus_{0\leq k\leq n}D_{[k+1,n+m]}\otimes
B_k}\leq\left( \bigoplus_{n+1\leq j\leq n+m}D_{[j+1,n+m]}\otimes
B_j\op \0\right) \ot S_n.
\end{eqnarray*}
Hence we have
\begin{eqnarray*}
{\bigoplus_{0\leq k\leq n+m}D_{[k+1,n+m]}\otimes
B_k}\leq \left( \bigoplus_{n+1\leq j\leq n+m}D_{[j+1,n+m]}\otimes
B_j\op \0\right) \ot S_n,
\end{eqnarray*}
from which the lemma follows since $S_{[n+1,n+m]}\geq 0$.
\endproof

\begin{lemma}\label{lem:lambda}
The following limit
\begin{eqnarray*}
\Lambda_S(\theta) =\lim_{n\to\infty}\frac{1}{n}\log \esp\left[
e^{\theta S_n}\right],
\end{eqnarray*}
exists in $\real\cup\{+\infty\}$ for all $\theta \geq 0$.
$\Lambda_S(.)$ is a proper convex function which is finite on the
interval $[0,\eta)$. Moreover for all $n$ and for all $\theta<\eta$,
we have $\esp[\exp \theta S_n]<\infty$.
\end{lemma}
\proof{Proof} In view of Lemma \ref{lem:aux}, we have clearly the
following subadditive property, for $\theta \geq 0$,
\begin{eqnarray*}
\log\esp\left[ e^{\theta S_{m+n}}\right]\leq \log\esp\left[e^{\theta
S_n} \right]+\log\esp\left[ e^{\theta S_{m-1}}\right],
\end{eqnarray*}
and the existence of $\Lambda_S$ follows and moreover, we have
\begin{eqnarray}\label{eq:inf}
\Lambda_S(\theta) = \inf_{n\geq 1}\frac{1}{n}\log\esp\left[e^{\theta
S_n} \right].
\end{eqnarray}
The last part of the lemma follows from
\begin{eqnarray*}
\bigoplus_i B_0^{(i)}\leq &S_n&\leq \sum_{k=0}^n \bigoplus_{i}
B^{(i)}_k,
\end{eqnarray*}
thanks to Assumption \ref{ass:sp}. Then we have
\begin{eqnarray*}
\left(\bigoplus_{i}\esp\left[
 e^{\theta B_0^{(i)}}\right]\right)\leq &\esp\left[e^{\theta
S_n}\right]&\leq \left(\bigotimes_{i}\esp\left[
 e^{\theta B_0^{(i)}}\right]\right)^{n+1},
\end{eqnarray*}
hence we see that $\Lambda_S(\theta)$ is finite for
$\theta\in[0,\eta)$ and infinite for $\theta>\eta$. Then the fact
that $\Lambda_S$ is a proper convex function follows from Lemma
2.3.9 of \cite{demzet}.
\endproof

\subsection{Computation of $\theta^*$}

\begin{lemma}\label{lem:theta}
We have $\theta^*>0$ and
\begin{eqnarray*}
\Lambda_S(\theta)+\Lambda_T(-\theta)<0 && \mbox{if}\quad \theta \in (0,\theta^*),\\
\Lambda_S(\theta)+\Lambda_T(-\theta)>0 && \mbox{if}\quad \theta
>\theta^*.
\end{eqnarray*}
\end{lemma}
\proof{Proof} In view of Assumption \ref{ass:s} and (\ref{eq:stab}),
we can choose $n$ such that $\esp[S_n]<na$. Then the function
$\theta\mapsto \log\esp[\exp(\theta S_n)]$ is convex, continuous and
differentiable on $[0,\eta)$. Hence we have
\begin{eqnarray*}
\frac{1}{n}\log\esp[e^{\theta S_n}]+\log\esp\left[e^{-\theta
\tau_1}\right] = \theta\left(\frac{\esp[S_n]}{n}-a\right)
+o(\theta),
\end{eqnarray*}
which is less than zero for sufficiently small $\theta>0$. By
(\ref{eq:inf}), we have for such $\theta$,
\begin{eqnarray*}
\Lambda_S(\theta) +\Lambda_T(-\theta)\leq
\frac{1}{n}\log\esp[e^{\theta S_n}]+\log\esp\left[e^{-\theta
\tau_1}\right]<0.
\end{eqnarray*}
Hence $\theta^*$ is positive and the last part of the lemma follows
from the convexity of $\Lambda_S(\theta)+\Lambda_T(-\theta)$.
\endproof

\subsection{Tail asymptotics for $Z$}

\begin{lemma}\label{lem:upper}
We have,
\begin{eqnarray*}
\limsup_{x\to \infty}\frac{1}{x}\log \prob(Z>x) \leq -\theta^*.
\end{eqnarray*}
\end{lemma}
\proof{Proof} First note that by Lemma \ref{lem:lambda}, we have
$\eta\geq \theta^*$ and for all $\theta<\eta$, we have
$\esp[\exp(\theta S_n)]<\infty$. We denote $S^\tau_n=(T_0-T_{-n})$
for $n\geq 0$. For any $0<\theta<\theta^*$, we have thanks to
Chernoff's inequality,
\begin{eqnarray*}
\prob(Z>x) &=&\prob \left(\cup_n\{S_n-S^\tau_n>x\}\right) \\
&\leq& \sum_n \prob(S_n-S^\tau_n>x)\\
&\leq& e^{-\theta x} \sum_n \esp\left[e^{\theta S_n}
\right]e^{n\Lambda_T(-\theta)},
\end{eqnarray*}
where the last series converge because we proved in Lemma
\ref{lem:theta} that for $\theta<\theta^*$, we have,
\begin{eqnarray*}
\frac{1}{n}\log\esp[e^{\theta S_n}]+\log\esp\left[e^{-\theta
\tau_1}\right]\rightarrow \Lambda_S(\theta) +\Lambda_T(-\theta)<0.
\end{eqnarray*}
\endproof

\begin{lemma}\label{lem:lower}
We have,
\begin{eqnarray*}
\liminf_{x\to \infty}\frac{1}{x}\log \prob(Z>x) \geq -\theta^*.
\end{eqnarray*}
\end{lemma}
\proof{Proof} We denote $\Lambda(\theta)=
\Lambda_S(\theta)+\Lambda_T(-\theta)$ and $S^\tau_n=(T_0-T_{-n})$ as
in previous proof. We only need to consider the case
$\theta^*<\infty$. We consider first the case where there exists
$\theta> \theta^*$ such that $\Lambda(\theta)<\infty$. First note
that the function $\theta\mapsto \Lambda(\theta)$ is convex, hence
the left-hand derivatives $\Lambda'(\theta-)$ and the right-hand
derivatives $\Lambda'(\theta+)$ exist for all $\theta>0$. Moreover,
we have $\Lambda'(\theta-)\leq\Lambda'(\theta+)$ and the function
$\theta\mapsto \frac{1}{2}(\Lambda'(\theta-)+\Lambda'(\theta+))$ is
non-decreasing, hence
$\Lambda'(\theta)=\Lambda'(\theta-)=\Lambda'(\theta+)$ except for
$\theta\in \Delta$, where $\Delta$ is at most countable. Since
$\Lambda(\theta)<\infty$ for $\theta>\theta^*$, we have
$\Lambda(\theta^*)=0$ and $\Lambda'(\theta^*+)>0$. To prove this,
assume that $\Lambda'(\theta^*+)=0$. Take $\theta<\theta^*$, thanks
to Lemma \ref{lem:theta}, we have $\Lambda(\theta)<0$. Choose
$\epsilon>0$ such that $0<\Lambda(\theta^*+\epsilon)<\epsilon
|\Lambda(\theta)|$. We have
\begin{eqnarray*}
\frac{\Lambda(\theta^*+\epsilon)}{\epsilon} <
\frac{-\Lambda(\theta)}{\theta^*-\theta},
\end{eqnarray*}
which contradicts the convexity of $\Lambda(\theta)$. Hence, we can
find $t\leq \theta^*+\epsilon$ such that
\begin{eqnarray*}
0<\Lambda(t),\quad t\notin \Delta.
\end{eqnarray*}
Note that these conditions imply $t>\theta^*$ and $\Lambda'(t)\geq
\Lambda'(\theta^*+)>0$.

Thanks to G\"artner-Ellis theorem (Theorem 2.3.6 in \cite{demzet}),
we have
\begin{eqnarray}
\label{eq:gart}\liminf_{n\rightarrow \infty} \frac{1}{n} \log
\prob(S_n-S^\tau_n>n \alpha) \geq - \inf_{x \in \mathcal{F},\:
x>\alpha} \Lambda^*(x),
\end{eqnarray}
where $\mathcal{F}$ is the set of exposed point of $\Lambda^*$ and
$\Lambda^*(x) = \sup_{\theta\geq 0}(\theta x-\Lambda(\theta))$. Note
that from the monotonicity of $\theta x-\Lambda(\theta)$ in $x$ as
$\theta$ is fixed, we deduce that $\Lambda^*$ is non-decreasing.
Moreover take $\alpha =\Lambda'(t)$, then $\Lambda^*(\alpha) =
t\alpha-\Lambda(t)$ and $\alpha\in \mathcal{F}$ by Lemma 2.3.9 of
\cite{demzet}.

Given $x>0$, define $n=\lceil x/\alpha\rceil$. We have
\begin{eqnarray*}
\frac{1}{x}\log \prob(Z>x) \geq
\frac{1}{n\alpha}\log\prob(S_n-S^\tau_n\geq n\alpha),
\end{eqnarray*}
taking the limit in $x$ and $n$ (while $\alpha=\Lambda'(t)$ is
fixed) gives thanks to (\ref{eq:gart}),
\begin{eqnarray*}
\liminf_{x\rightarrow \infty} \frac{1}{x} \log
\prob(Z>x)\geq-\frac{t\alpha-\Lambda(t)}{\alpha} \geq -t\geq
-\theta^*-\epsilon.
\end{eqnarray*}

We consider now the case where for all $\theta>\theta^*$, we have
$\Lambda(\theta)=\infty$, hence $\eta=\theta^*$. Take $K>0$ and define $\tilde{S}^K_{[n,m]} =
S_{[n,m]}\prod_{i=n}^m\ind(\oplus_j B_i^{(j)}\leq K)$ and $\tilde{Z}^K =
\sup_{n\geq 0} (\tilde{S}^K_{[-n,0]}-S_n^\tau)$. We have clearly $Z\geq
\tilde{Z}^K$. It is easy to see that the proof of Lemma
\ref{lem:lambda} is still valid (note that the subadditive property
carries over to $\tilde{S}^K_{[n,m]}$) and the following limit exists
\begin{eqnarray*}
\tilde{\Lambda}^K_S(\theta) = \lim_{n\to \infty}\frac{1}{n}\log\esp\left[
  e^{\theta \tilde{S}^K_{[1,n]}}\right] = \inf_n\frac{1}{n}\log\esp\left[
  e^{\theta \tilde{S}^K_{[1,n]}}\right].
\end{eqnarray*}
Moreover thanks to the subadditive property of $S$ (see Lemma
\ref{lem:aux}), we have $\tilde{S}^K_{[1,n]}\leq
\tilde{S}^K_{[1,1]}+\dots+\tilde{S}^K_{[n,n]}=\oplus_j
B_1^{(j)}+\dots+\oplus_j B_n^{(j)}$. Hence we have
$\prob(\tilde{S}^K_{[1,n]}\leq nK)=1$, so that
$\tilde{\Lambda}^K_S(\theta)\leq \theta K$. Hence by the first part of
the proof, we have
\begin{eqnarray*}
\liminf_{x\to\infty}\frac{1}{x}\log\prob(\tilde{Z}^K>x) \geq -\tilde{\theta}^K,
\end{eqnarray*}
with $\tilde{\theta}^K = \sup\{\theta>0,\:
\tilde{\Lambda}^K_S(\theta)+\Lambda_T(-\theta)<0\}$. We now prove that
$\tilde{\theta}^K\to\eta$ as $K$ tends to infinity which will
conclude the proof.
Note that for any fixed $\theta\geq 0$, the function
$\tilde{\Lambda}^K_S(\theta)$ is nondecreasing in $K$ and
$\lim_{K\rightarrow \infty}\tilde{\Lambda}^K_S(\theta) =
\tilde{\Lambda}_S(\theta)\leq \Lambda_S(\theta)$. This directly implies
that $\tilde{\theta}^K\geq \eta$.
Take $\theta>\eta$, so that $\Lambda_S(\theta)=\infty$. If
$\tilde{\Lambda}_S(\theta)<\infty$, then for all $K$, we have
$\tilde{\Lambda}^K_S(\theta)\leq\tilde{\Lambda}_S(\theta)<\infty$. But, we
have $\tilde{\Lambda}^K_S(\theta)=\inf_n
\frac{1}{n}\log\esp\left[
  e^{\theta \tilde{S}^K_{[1,n]}}\right]$, so that there exists $n$ such that
\begin{eqnarray*}
\esp\left[ e^{\theta S_{[1,n]}},\:\max(\oplus_j B^{(j)}_1,\dots ,\oplus_jB^{(j)}_n)\leq K\right]\leq e^{\tilde{\Lambda}^K_S(\theta)+1}\leq e^{\tilde{\Lambda}_S(\theta)+1},
\end{eqnarray*}
but the left-hand side tends to infinity as $K\to \infty$. Hence we
proved that for all $\theta>\eta$, we have
$\tilde{\Lambda}^K_S(\theta)\to\infty$ as $K\to \infty$. This implies
that $\tilde{\theta}^K\to \eta$ as $K\to \infty$.
\endproof

\section{Proof of Theorem \ref{the:2} and Corollary \ref{cor}}\label{sec:p2}

We begin with a general result showing the existence of the function
$\Lambda_\ell$. Let $\{M_n\}$ be an i.i.d. sequence of irreducible
aperiodic (max,plus)-matrices with fixed structure. We denote
\begin{eqnarray*}
M_{[1,n]}^{(i,j)} &=& \left( M_n\otimes\dots \otimes
M_1\right)^{(i,j)}.
\end{eqnarray*}

\begin{lemma}\label{lem:M}
For $\theta\geq 0$, the following limit exists in
$\real\cup\{+\infty\}$ and is independent of $i$ and $j$,
\begin{eqnarray*}
\Lambda_M(\theta) = \lim_{n\rightarrow
\infty}\frac{1}{n}\log\esp\left[e^{\theta M_{[1,n]}^{(i,j)}}\right].
\end{eqnarray*}
\end{lemma}
\proof{Proof}

We denote
\begin{eqnarray*}
\Lambda_M^{(i,j)}(\theta,n) &=& \log\esp\left[e^{\theta
M_{[1,n]}^{(i,j)}}\right].
\end{eqnarray*}
We have for $\theta\geq 0$,
\begin{eqnarray*}
\Lambda^{(i,j)}_M(\theta,n+m) &=& \log\esp\left[e^{\theta M_{[1,n+m]}^{(i,j)}}\right]\\
&=& \log\esp\left[\max_k e^{\theta M_{[n+1,n+m]}^{(i,k)}}e^{\theta M_{[1,n]}^{(k,j)}}\right]\\
&\geq& \max_k\left\{\log\esp\left[e^{\theta M_{[n+1,n+m]}^{(i,k)}}\right]+\log\esp\left[e^{\theta M_{[1,n]}^{(k,j)}}\right]\right\}\\
&=&
\max_k\left\{\Lambda_M^{(i,k)}(\theta,m)+\Lambda_M^{(k,j)}(\theta,n)\right\}.
\end{eqnarray*}
In particular for $j=i$, we have
\begin{eqnarray*}
\Lambda_M^{(i,i)}(\theta,n+m) &\geq&
\Lambda_M^{(i,i)}(\theta,m)+\Lambda_M^{(i,i)}(\theta,n).
\end{eqnarray*}
Moreover thanks to the fixed structure assumption and the
aperiodicity, there exists $N$ such that for $n\geq N$, we have
$M_{[1,n]}^{(i,j)}>-\infty$ for all $i$ and $j$, hence
$\Lambda_M^{(i,j)}(\theta,n)>-\infty$ and we have
\begin{eqnarray*}
\lim_{n\rightarrow \infty}\frac{1}{n}
\Lambda_M^{(i,i)}(\theta,n)=\sup_{n\geq
N}\frac{1}{n}\Lambda_M^{(i,i)}(\theta,n)>-\infty.
\end{eqnarray*}
For arbitrary $i$ and $j$, choose $n,m\geq N$ and note that
\begin{eqnarray*}
\Lambda_M^{(i,j)}(\theta,n+m)&\geq &\Lambda_M^{(i,i)}(\theta,n)+\Lambda_M^{(i,j)}(\theta,m),\\
\Lambda_M^{(i,i)}(\theta,n+m)&\geq
&\Lambda_M^{(i,j)}(\theta,n)+\Lambda_M^{(j,i)}(\theta,m),
\end{eqnarray*}
where all terms are in $\real\cup\{+\infty\}$. Letting $n\rightarrow
\infty$ while keeping $m$ fixed, it follows that
\begin{eqnarray*}
\lim_{n\rightarrow\infty}
\frac{1}{n}\Lambda_M^{(i,j)}(\theta,n)=\lim_{n\rightarrow\infty}
\frac{1}{n}\Lambda_M^{(i,i)}(\theta,n).
\end{eqnarray*}
\endproof

Note that Assumption {\bf (ST)} ensures that the matrices
$\{A_n(\ell,\ell)\}$ are irreducible aperiodic with fixed structure.
We now extend previous lemma to the sequence $\{A_n\}$ of reducible
matrices. To do so, we first associate a graph $\Gcal=(\Vcal,\Ecal)$
to $A_n$, as in Section~2.3 of \cite{bcoq}. Set
$\Vcal:=\{1,\dots,s\}$, which we abbreviate as $[1,s]$. An edge
$(i,j)$ belongs to $\Ecal$ if and only if $A_n^{(j,i)}\geq 0$. Two
nodes of $\Vcal$ are said to belong to the same communication class
if there is a directed path from the first to the second and another
one from the second to the first. Let $\Ccal_1,\dots, \Ccal_{d}$ be
the communication classes of $\Gcal$ and $\lessdot$ the associated
partial order, namely $\Ccal_\ell\lessdot \Ccal_m$ if there is a
path from any vertex in $\Ccal_\ell$ to any vertex in $\Ccal_m$.
Without loss of generality, we assume that $\mathcal{C}_\ell\lessdot
\mathcal{C}_m$ implies $\ell\leq m$; this is a notationally
convenient restriction on the numbering of the communication
classes.

We use the following notation:
\begin{itemize}
\item for any coordinate $i\in \mathcal{V}$, its communication class is denoted by
$[i]$,
\item for any coordinate $i$,
the subset of coordinates $j$ such that $[j]\lessdot [i]$ is denoted
by $[\leq i]$;
\item for any coordinate $i$,
the subset of coordinates $j$ such that $[i]\lessdot [j]$ is denoted
by $[i\leq]$;
\item for any coordinate $i$ and $j\in[i\leq]$, we
write
$$
[i\leq j]:= [i\leq]\cap [\leq j].
$$
\end{itemize}

We now extend previous lemma. We introduce first some notations,
\begin{eqnarray*}
\Lambda_{\ell}(\theta)=\Lambda_{[i]}(\theta)&=&\lim_{n\rightarrow
\infty}\frac{1}{n}\log\esp \left[ e^{\theta
D_{[1,n]}^{(i,i)}}\right],
\end{eqnarray*}
which does not depend on $i\in \mathcal{C}_\ell$ as shown above.
\begin{lemma}
For $\theta\in [0,\eta)$, we have
\begin{eqnarray*}
\lim_{n\to\infty}\frac{1}{n}\log\esp\left[ e^{\theta
D_{[1,n]}^{(i,j)}}\right] = \sup_{k \in [i\leq
j]}\Lambda_{[k]}(\theta),
\end{eqnarray*}
where the supremum over the empty set is $-\infty$.
\end{lemma}
\proof{Proof} If $[i\leq j] =\emptyset$, the result is obvious and
if $j\in [i]$, the result follows from previous lemma. Hence we
consider only the case: $[i\leq j]\supset [i]$. We denote
\begin{eqnarray*}
\Lambda^{(i,j)}(\theta,n) &=& \log\esp\left[e^{\theta
D_{[1,n]}^{(i,j)}}\right].
\end{eqnarray*}
With the same argument as in previous lemma, we have for $\theta\geq
0$, and for any $[i]\lessdot \Ccal_\ell\lessdot [j]$, there exists
$k\in \Ccal_\ell$ and $u,v\geq 1$,
\begin{eqnarray*}
{\Lambda^{(i,j)}(\theta,n+u+v) \geq}\Lambda^{(i,k)}(\theta,u)+\Lambda^{(k,k)}(\theta,n)+\Lambda^{(k,j)}(\theta,v),
\end{eqnarray*}
where each term is finite since $\theta<\eta$. Hence we have by
previous lemma:
\begin{eqnarray*}
\liminf_{n\to\infty}\frac{\Lambda^{(i,j)}(\theta,n)}{n}\geq\sup_{k
\in [i\leq j]}\Lambda_{[k]}(\theta).
\end{eqnarray*}
For the upper bound, note that there exists $u\geq 1$ such that
$D^{(i,j)}_{[1,u]}>-\infty$ for all $i\leq j$. Consider first the
case $d=2$, $i\in \Ccal_1$ and $j\in \Ccal_2$, then we have
\begin{eqnarray*}
\Lambda^{(i,j)}(\theta,n) &\leq &\log\esp\left[\max_{a\in\Ccal_2}
e^{\theta D^{(i,a)}_{[1,u]}}e^{\theta
D^{(a,j)}_{[u+1,n]}}+\max_{b\in\Ccal_1} e^{\theta D^{(i,b)}_{[1,n-u]}}e^{\theta
D^{(b,j)}_{[n-u+1,n]}}\right],\\
&\leq& \log\left(\sum_{a\in\Ccal_2}e^{\Lambda^{(i,a)}(\theta,u)}
e^{\Lambda^{(a,j)}(\theta,n)}+\sum_{b\in\Ccal_1}e^{\Lambda^{(i,b)}(\theta,n)}
e^{\Lambda^{(b,j)}(\theta,u)}\right),
\end{eqnarray*}
Hence by Lemma 1.2.15 of \cite{demzet}, we have
\begin{eqnarray*}
\limsup_{n\to \infty}\frac{\Lambda^{(i,j)}(\theta,n)}{n} \leq
\max\left( \Lambda_{1}(\theta),\Lambda_2(\theta)\right).
\end{eqnarray*}
We have clearly by induction that
\begin{eqnarray*}
\limsup_{n\to \infty}\frac{\Lambda^{(i,j)}(\theta,n)}{n} \leq
\sup_{k\in[i\leq j]} \Lambda_{[k]}(\theta),
\end{eqnarray*}
which concludes the proof.
\endproof

We now compute $\Lambda_S(\theta)$ for a (max,plus)-linear system
under the assumptions of Theorem \ref{the:2}.
\begin{lemma}\label{lem:lambdamaxplus}
We have for $\theta \in [0,\eta)$
\begin{eqnarray*}
\Lambda_S(\theta) = \sup_\ell \Lambda_\ell(\theta).
\end{eqnarray*}
\end{lemma}
\proof{Proof}

The lower bound follows directly from the following inequality: for
all $\ell$, we have for $i\in \Ccal_\ell$,
\begin{eqnarray*}
\esp\left[ e^{\theta S_n}\right]&\geq&\esp\left[ e^{\theta
D_{[1,n]}^{(i,i)}}\right].
\end{eqnarray*}
We now derive the upper bound. Note that $A_k\ot \0 =B_k\op \0\geq
B_k$, hence we have
\begin{eqnarray*}
S_n&=&\bigoplus_{1\leq i\leq s}\bigoplus_{0\leq k\leq
n}\left(D_{[k+1,n]}\otimes
B_k\right)^{(i)}\\
&\leq&\bigoplus_{1\leq i,j\leq s}\bigoplus_{0\leq k\leq
n}D_{[k,n]}^{(i,j)}.
\end{eqnarray*}
Hence we have
\begin{eqnarray*}
\esp\left[ e^{\theta S_n}\right]&\leq& \sum_{i,j}\sum_k \esp\left[
e^{\theta D_{[k,n]}^{(i,j)}}\right]\\
&\leq& n\sum_{i,j} \esp\left[ e^{\theta D_{[0,n]}^{(i,j)}}\right],
\end{eqnarray*}
and the lemma follows directly from Lemma 1.2.15 of \cite{demzet}.
\endproof

Theorem \ref{the:2} follows directly form the fact that
$\Lambda_S(\theta)=\infty$ as soon as $\theta>\eta$ which follows
from the lower bound $B_0^{(i)}\leq S_n$.

We now prove Corollary \ref{cor}. The following lemma implies that
$\Lambda_S(\theta) = \sup_\ell \Lambda_\ell(\theta)$ for all
$\theta\geq 0$ and then Corollary \ref{cor} follows:
\begin{lemma}
Under assumptions of Corollary \ref{cor}, there exists $\ell\in
[1,d]$ such that $\Lambda_\ell(\theta) = \infty$ for all
$\theta>\eta$.
\end{lemma}
\proof{Proof} We only need to consider the case $\eta<\infty$. Take
$i$ such that $\sup\{\theta,\:\esp[\exp\theta
B^{(i)}_0]<\infty\}=\eta$. By the condition on the entries of $A_n$,
we have
\begin{eqnarray*}
B_n^{(i)}\leq \oplus_jA_n^{(i,j)} \leq \sum_{k=1}^K \sigma^{(k)}_n.
\end{eqnarray*}
Hence there exists $k$ such that $\esp[\exp(\theta
\sigma^{(k)}_n)]=\infty$ for $\theta>\eta$ and there exist $j$ such
that $A_n^{(j,j)}=\sigma^{(k)}_n$. Then $\ell$ defined by
$\Ccal_\ell=[j]$ satisfies the property claimed in the lemma.
\endproof

\section{Conclusion}

We have shown that the distribution of the stationary solution of a
(max,plus) recursion has an exponentially decaying tail and we gave
an analytical way to compute the decay rate.

We applied our results to different kinds of communication networks
and exhibited quite non-standard behavior possible in high-dimension
only. We also analyzed a queueing network with multipath routing and
showed on a simple example how our analysis could help in the design
of the routing decision depending on the characteristic of the
traffic.

We are currently working on some possible extensions of our
work. Of particular interest would be a large deviations principle for
the process $\{S_n/n\}$ introduced in (\ref{eq:defS_n}). It would
allow to give the most probable way for a large deviation of the
maximal dater.
Also, it
should be possible to use the distributional Little's law \cite{hg}
to get asymptotics for the number of packets in the networks.

In general, the characterization of the decay rate is given by the
moment generating function which is not easy to compute, especially
in the case of feedback. One practical question of interest would be
to find good ways to estimate this function from the statistics made
on the traffic.

\section*{Acknowledgment}

The author would like to thank the participants of Valuetools 2006
(where this work was presented) and especially Bruno Gaujal for a
comment related to Remark \ref{rem:bruno}. He would also like to
thank Peter Friz for pointing out a mistake in an earlier version of
this paper.

\bibliography{ex}
\bibliographystyle{chicago}

Marc Lelarge\footnote{This work was partially done while the author
  was with Boole Centre for Research in Informatics, Science
Foundation Ireland Research Grant No. SFI 04/RP1/I512.}\\
ENS-INRIA\\
45 rue d'Ulm\\
75005 Paris, France\\
\\  {\tt e-mail : marc.lelarge@ens.fr}\\
 {\tt http://www.di.ens.fr/$\sim$lelarge}

\end{document}